\documentclass[a4paper,10pt]{amsart}
\usepackage[utf8]{inputenc}
\usepackage{lmodern}
\usepackage[english]{babel}
\usepackage[T1]{fontenc}
\usepackage{microtype} 
    \usepackage[OT2,OT1]{fontenc}
    \newcommand\cyr{%
    \renewcommand\rmdefault{wncyr}%
    \renewcommand\sfdefault{wncyss}%
    \renewcommand\encodingdefault{OT2}%
    \normalfont
    \selectfont}
    \DeclareTextFontCommand{\textcyr}{\cyr}
\usepackage{amsmath,amssymb,amsfonts,amsthm,amscd,latexsym,mathrsfs}
\usepackage[all]{xy}
\input xypic
\usepackage{color}
\usepackage{hyperref}
\usepackage{lipsum}
\usepackage{breqn}
\hypersetup{colorlinks=true,linkcolor=red,citecolor=blue}
\usepackage{booktabs}
\usepackage{multirow}
\theoremstyle{plain}
\newtheorem{theorem}[subsection]{{\bf Theorem}}
\newtheorem*{theorem*}{{\bf Theorem}}
\newtheorem{corollary}[subsection]{{\bf Corollary}}
\newtheorem*{corollary*}{{\bf Corollary}}
\newtheorem{proposition}[subsection]{{\bf Proposition}}

\newtheorem{question}[subsection]{{\bf Question}}

\theoremstyle{definition}

\theoremstyle{remark}

\numberwithin{equation}{section}
\DeclareMathOperator{\Ima}{Im}

\DeclareMathOperator{\kernal}{ker}
\newcommand{\gen}[1]{\langle #1 \rangle}

\begin{document}
\baselineskip=14pt
\title{Bounds On the order of the Schur multiplier of $p$-groups}
\author[P. K. Rai]{Pradeep K. Rai}
\address[Pradeep K. Rai]{Mahindra University, Hyderabad, Telangana,\\
India}
\email{raipradeepiitb@gmail.com}
\begin{abstract}

In 1956, Green provided a bound on the order of the Schur multiplier of $p$-groups. This bound, given as a function of the order of the group, is the best possible. Since then, the bound has been refined numerous times by adding other inputs to the function, such as, the minimal number of generators of the group and the order of the derived subgroup. We strengthen these bounds by adding another input, the group's nilpotency class. The specific cases of nilpotency class 2 and maximal class are discussed in greater detail.

\end{abstract}
\maketitle
\section{Introduction}





\vspace{.2cm}
The Schur multiplier $M(G)$ of a finite group $G$ is defined as the second cohomology group of $G$ with coefficients in $\mathbb{C}^{*}$. It plays an important role in the theory of extensions of groups. Finding the bounds on the order, exponents, and ranks of the Schur multiplier of prime power groups has been a major focus of previous investigations. This article investigates the bounds on the order of the Schur multiplier of prime power groups.

Let $G$ be a finite $p$-group of order $p^n$. In 1956 Green \cite{Green} proved that $|M(G)| \leq p^{\frac{1}{2}n(n-1)}$. Since then, this bound has been strengthened  by many mathematicians \cite{Blackburn,Ellis1,Gaschutz,Jones2,Jones, Jones1,Moravec,Niroomand2,Niroomand1,Rai1,Vermani1,Vermani2,Weigold1,Weigold2}. 

To note the most recent ones, let G be a non-abelian $p$-group of order $p^n$ with derived subgroup of order $p^k$.  Niroomand \cite{Niroomand2} proved that 
\begin{equation} \label{eq0}
|M(G)| \leq p^{\frac{1}{2}(n-k-1)(n+k-2)+1}. 
\end{equation}

The author noted in \cite{Rai1} that a bound by Ellis and Wiegold \cite{Ellis1} is better than this bound and derived from their bound that 
\[|M(G)| \leq p^{\frac{1}{2}(d-1)(n+k-2)+1}(= p^{\frac{1}{2}(d-1)(n+k)-(d-2)}),\]  
where $d$ is the minimal number of generators of $G.$

In this article, we further refine the bounds by adding another input, i.e., the nilpotency class of the group. 

Before proceeding to the results of this article, we set some notations that are mostly standard. The center and the commutator subgroup of a group $G$ are denoted by $Z(G)$ and $\gamma_2(G)$, respectively. By $d(G)$ we denote the minimal number of generators of $G$. We write $\gamma_i(G)$ for the $i$-th term in the lower central series of $G$. The subgroup $\gen{x^{p} \mid x \in G}$ is denoted by $G^{p}$. Finally, the abelianization of the group $G$, i.e. $G/\gamma_2(G)$, is denoted by $G^{ab}$. \\

We now state our first theorem.
\begin{theorem}\label{thm1}
Let $G$ be a non-abelian $p$-group of order $p^n$ and nilpotency class $c$ with $|\gamma_2(G)| = p^k$ and $d =d(G)$. Then
\[|M(G)| \leq p^{\frac{1}{2}(d-1)(n+k)-\sum\limits_{i =2}^{\min(d,c)} d-i}.\] 
Thus, if $\mu = \min(d, c)$ then, 
\[|M(G)| \leq p^{\frac{1}{2}(d-1)(n+k)-\frac{1}{2}(\mu-1)[2d-(\mu+2)]}.\]
Considering cases, the inequality can be restated as follows:
\[ |M(G)| \leq \begin{cases} 
        p^{\frac{1}{2}(d-1)(n+k) - \frac{1}{2}(d-1)(d-2)} & \text{if} \ \ \ \ \  d \leq c  \\ \\
      p^{\frac{1}{2}(d-1)(n+k)- \frac{1}{2}(c-1)(2d-(c+2))} & \text{if} \ \ \ \ \  d \geq c. \\
   \end{cases}
\]
\end{theorem}

\vspace{.6cm}

Next, we consider $p$-groups of nilpotency class 2. A finite $p$-group of nilpotency class 2 is said to be special if its center coincides with the derived and the Frattini subgroups. Berkovich and Janko asked the following questions:

\begin{question}\cite[Problem 1729]{Berkovich}\label{q1}
 Let G be a special p-group with $d(G) = d$ and $|Z(G)| = p^{\frac{1}{2}d(d-1)}$. Find the Schur multiplier of $G$ and describe the representation groups of $G$.
\end{question}

\begin{question}\cite[Problem 2027]{Berkovich}\label{q2}
 Find the Schur multiplier of special $p$-groups with center of order $p^2$.
\end{question}
The Question \ref{q1} and \ref{q2} of Bekovich and Janko have been studied in \cite{Rai2} and \cite{Hatui} respectively. 
In view of these questions, and the fact that for special $p$-groups $G$, $|Z(G)| = p^k$ if and only if $d(\gamma_2(G)) = k$, it seems reasonable to consider the term $d(\gamma_2(G))$ while investigating the bounds for the order of the Schur multiplier of $p$-groups of nilpotency class 2. 

\begin{theorem}\label{thm2}
 Let $G$ be a non-abelian finite $p$-group of order $p^n$ with $|\gamma_2(G)| =p^k$, $d(G) = d$ 
and $d\bigg(\frac{\gamma_2(G)}{\gamma_3(G)}\bigg) = \gamma$. 
 Then
\[|M(G)| \leq p^{\frac{1}{2}(d-1)(n+k)-\sum\limits_{i =2}^{\min(d,\gamma+1)} d-i}.\]
Thus, if $\nu = \min(d, \gamma+1)$ then, 
\[|M(G)| \leq p^{\frac{1}{2}(d-1)(n+k)-\frac{1}{2}(\nu-1)[2d-(\nu+2)]}.\]
Considering cases, the inequality can be restated as follows:
\[ |M(G)| \leq \begin{cases} 
        p^{\frac{1}{2}(d-1)(n+k) - \frac{1}{2}(d-1)(d-2)} & \text{if} \ \ \ \ \  d \leq \gamma+1  \\ \\
      p^{\frac{1}{2}(d-1)(n+k)- \frac{1}{2}\gamma(2d-(\gamma+3))} & \text{if} \ \ \ \ \  d \geq \gamma+1. \\
   \end{cases}
\]
\end{theorem}  

 

In view of the Questions \ref{q1} and \ref{q2} the following corollary is an application of Theorem \ref{thm2} to special $p$-groups.
 \begin{corollary}\label{cor2}
 Let $G$ be a finite $p$-group of nilpotency class 2 such that $G^p \leq \gamma_2(G)$, $|\gamma_2(G)| = p^{k}$ and $d(G) = d$. Then 
 \[ p^{\frac{1}{2}d(d-1) - k} \leq |M(G)| \leq \begin{cases} 
        p^{(d-1)(k+1)} & \text{if} \ \ \ \ \  d \leq k+1  \\ \\
      p^{\frac{1}{2}d(d-1)+ \frac{1}{2}k(k+1)} & \text{if} \ \ \ \ \  d \geq k+1.  \\
   \end{cases}
\]
Moreover 
\begin{itemize}
    \item 
If $G^p \cong \mathbb{Z}_p$ then, 
 \[ |M(G)| \leq \begin{cases} 
        p^{(d-1)k-1} & \text{if} \ \ \ \ \  d \leq k  \\ \\
      p^{\frac{1}{2}d(d-1)+ \frac{1}{2}k(k-1)-1} & \text{if} \ \ \ \ \  d \geq k. \\
   \end{cases}
\]
\item 
If $p$ is an odd prime and $G^p = \gamma_2(G)$, then
\[|M(G)| \leq p^{\frac{1}{2}d(d-1)+ \frac{1}{2}k(k-3)}.\]
\item 
If $G_1$ and $G_2$ are special $p$-groups with $|Z(G_1)| = p^2$ and $|Z(G_2)| = p^3$ then
\[p^{\frac{1}{2}d(d-1) - 2} \leq |M(G_1)| \leq p^{\frac{1}{2}d(d-1) +3}\]
and 
\[p^{\frac{1}{2}d(d-1) - 3} \leq |M(G_2)| \leq p^{\frac{1}{2}d(d-1) +6}.\]
\end{itemize}
\end{corollary}

 Let $p$ be an odd prime and $G_1$ be a special $p$-group with $|Z(G_1)| = p^2$. A more general result has already been given by Mazur \cite{Mazur} in this case. He proves that if the epicenter $Z^{*}(G_1)$ does not coincide with the center $Z(G_1)$ then $G_1$ is of exponent $p$ and belongs to one of the five classes of groups. These classes include a group of order $p^5$, three groups of order $p^6$, two groups of order $p^7$, a group of order $p^{2m+3}$ (for all $m \geq 3)$ and two groups of order $p^{2m+2}$ (for all $m \geq 2$). 
It is easy to see from \cite[Theorem 2.5.10]{Karp} that if $Z^{*}(G_1)$ coincide with $Z(G_1)$ then $M(G_1)$ is elementary abelian of order $p^{\frac{1}{2}d(d-1) - 2}$. Therefore, it only remains to compute the Schur multiplier of groups that belong to the above mentioned five families. This can be easily achieved by using the explicit presentation of the non-unicentral groups given by Mazur. This observation nullifies \cite[Theroem 1.3(d)]{Hatui} which states that if $G_1^p \cong \mathbb{Z}_p$ then $|M(G_1)| = p^{\frac{1}{2}d(d-1)}$ if and only if $Z^{*}(G_1) = G^p$. It is clear from the above discussion that such groups do not exist. Also, in contrast to \cite[Theorem 1.1(c)]{Hatui}, one can see that the Schur multiplier of $G_1$ is always elementary abelian and never of exponent $p^2$.\\

The following corollary improves on Corollary \ref{cor2} for $G_2$ when $|G_2| \geq p^{13}.$

\begin{corollary}\label{cor3}
Let $p$ be an odd prime and $G_2$ be a special $p$-group with $|Z(G_2)| = p^3$ and $|G_2| \geq p^{13}.$ Then 
\[|M(G_2)| \leq p^{\frac{1}{2}d(d-1)+ 2}.\]
Moreover, if $G_2^p = \gamma_2(G_2)$, then
\[|M(G_2)| \leq p^{\frac{1}{2}d(d-1) - 2}.\]

\end{corollary}

Next, we consider the groups of maximal class. A finite $p$-group of order $p^n$ is said to be of maximal class if its nilpotency class is $n-1$. Let $G$ be a finite $p$-group of maximal class and order $p^n$. Since $G$ is generated by 2 elements, it follows by a result of Gasch\"{u}tz \cite{Gaschutz} that $|M(G)| \leq p^{n-1}$. Moravec \cite{Moravec} proved that, if $n > p+1$ then $|M(G)| \leq p^{\frac{p+1}{2}\left \lceil{\frac{n-1}{p-1}}\right \rceil}$. Improving Moravec's result, we prove the following theorem.

\begin{theorem} \label{thm3}
Let $p$ be an odd prime and $G$ be a finite $p$-group of maximal class with $ |G| = p^n$, $n \geq 4$.  Then $|M(G)| \leq p^{\frac{n}{2}}$.  
\end{theorem}

\section{Proofs}
Let $G$ be a finite $p$-group and $\overline{G}$ be the factor group $G/Z(G)$. The commutator $x^{-1}y^{-1}xy$ of the elements $x, y \in G$ is denoted by $[x,y].$ For ease of reading, we shall use the same `bar notation' $\overline{g}$ for $g \in G$ to denote the different elements in different factor groups, when there is no danger of ambiguity. Such notations should be interpreted according to the context. For example, whenever $\overline{[x_1, x_2]} \otimes \overline{x_{3}} \in  \frac{\gamma_2(G)}{\gamma_{3}(G)} \otimes \overline{G}^{ab}$ for $x_1, x_2, x_3 \in G$, by $\overline{[x_1, x_2]}$ and $\overline{x_{3}}$ we mean $[x_1,x_2]\gamma_3(G)$ and $x_3Z(G)\gamma_2(G)$, respectively.\\

We now proceed to prove the Theorem \ref{thm1}. The proof is founded on the following result of Ellis and Weigold \cite[ Proposition 1, comments following Theorem 2]{Ellis1}.

\begin{proposition}\label{prop1}
Let $G$ be a finite $p$-group of nilpotency class $c$ and $\overline{G}$ be the factor group $G/Z(G)$. Then 
\[\Big{|}M(G)\Big{|}\Big{|}\gamma_2(G)\Big{|}\prod_{i=2}^{c}\Big{|}\Ima \Psi_i\Big{|} \leq \Big{|}M(G^{ab})\Big{|}\prod_{i=2}^{c}\Big{|}\frac{\gamma_i(G)}{\gamma_{i+1}(G)} \otimes \overline{G}^{ab}\Big{|},\]

where $\Psi_i$, for $i=2, \ldots c,$ is a map from $\underbrace{\overline{G}^{ab} \otimes \overline{G}^{ab} \cdots \otimes \overline{G}^{ab}}_{i+1 \ \text{times}}$ to $\frac{\gamma_i(G)}{\gamma_{i+1}(G)} \otimes \overline{G}^{ab}$  defined as follows:\\
\[ \Psi_2(\overline{x_1} \otimes \overline{x_2} \otimes \overline{x_{3}}) 
= \overline{[x_1, x_2]} \otimes \overline{x_{3}} + \overline{[x_2, x_{3}]} \otimes \overline{x_1} + \overline{[x_3, x_{1}]} \otimes \overline{x_2}.\]
For $3 \leq i \leq c$,
\begin{eqnarray*}
\Psi_i(\overline{x_1} \otimes \overline{x_2} \otimes \cdots \otimes \overline{x_{i+1}}) 
& = & \overline{[x_1, x_2, \cdots, x_i]_l} \otimes \overline{x_{i+1}} + \overline{[x_{i+1}, [x_1, x_2, \cdots x_{i-1}]_l]} \otimes \overline{x_i} \\
&& +\overline{[[x_i, x_{i+1}]_r, [x_1, \cdots, x_{i-2}]_l]} \otimes \overline{x_{i-1}} \\
&& + \overline{[[x_{i-1}, x_i, x_{i+1}]_r, [x_1, x_2, \cdots, x_{i-3}]_l]} \otimes \overline{x_{i-2}} \\
&& + \cdots + \overline{[x_2, \cdots, x_{i+1}]_r} \otimes \overline{x_1}\\
\end{eqnarray*}
where 
\[[x_1, x_2, \cdots x_i]_r = [x_1, [\cdots [x_{i-2},[x_{i-1},x_i]]\ldots]\]
and 
\[[x_1, x_2, \cdots x_i]_l = [\ldots[[x_1, x_2], x_3], \cdots, x_i].\]
\end{proposition}

We are now ready to prove Theorem \ref{thm1}.\\

\noindent \textbf{\textit{Proof of Theorem \ref{thm1}}}:
Let $\Psi_i$ be the map defined above and $d(G/Z(G)) = \delta$. Following Proposition \ref{prop1} we have that 
\[|M(G)||\gamma_2(G)|\prod_{i=2}^{c}|\Ima \Psi_i| \leq |M(G^{ab})|p^{k\delta}.\]
Applying \cite[Lemma 2.1]{Rai1} this gives
\[|M(G)|\prod_{i=2}^{c}|\Ima \Psi_i| \leq p^{\frac{1}{2}(d-1)(n-k)+k(\delta-1)},\]
so that
\begin{equation}\label{eq1}
|M(G)|\prod_{i=2}^{c}|\Ima \Psi_i| \leq p^{\frac{1}{2}(d-1)(n+k) - k(d-\delta)}.
\end{equation}
Choose a subset $S= \{x_1, x_2, \ldots, x_{\delta}\}$ of $G$ such that $\{\overline{x_1},\overline{x_2}, \cdots, \overline{x_{\delta}}\}$ be a minimal generating set for $G/Z(G)$. Fix $i \leq \min(\delta, c)$. Since $ i \leq c$, $\gamma_i(G)/\gamma_{i+1}(G)$ is a non-trivial group. Using \cite[Lemma 3.6 (c)]{Khukhro} we can choose a commutator $[y_1,y_2, \cdots, y_i]$ of weight $i$ such that $[y_1,y_2, \cdots, y_i] \notin \gamma_{i+1}(G)$ and $y_1, \ldots, y_i \in S$ . Since $i \leq \delta$, $S \backslash \{y_1,y_2, \cdots, y_i\}$ contains at least $\delta - i$ elements. Choose any $\delta - i$ elements $z_1, z_2, \ldots, z_{\delta - i}$ from $S \backslash \{y_1,y_2, \cdots, y_i\}$. Since $[y_1,y_2, \cdots, y_i] \notin \gamma_{i+1}(G)$ and $z_j \notin \{y_1,y_2, \cdots, y_i\}$, $\Psi_i(\overline{y_1}, \ldots, \overline{y_i}, \overline{z_j}) \neq 1$. Notice that the set $\{\Psi_i(\overline{y_1}, \ldots, \overline{y_i}, \overline{z_j}) \ \ | \ \ 1 \leq j \leq \delta -i\}$ is a minimal generating set for $\gen{\{\Psi_i(\overline{y_1}, \ldots, \overline{y_i}, \overline{z_j}) \ \ | \ \ 1 \leq j \leq \delta -i\}}$, because $\{\overline{x_1},\overline{x_2}, \cdots, \overline{x_{\delta}}\}$ is a minimal generating set for $G/Z(G)$. It follows that $|\text{Im}\Psi_i| \geq p^{\delta - i}$. Putting this in Equation \ref{eq1} we get the required result.

 \vspace{.3cm}
     
We now proceed to prove Theorem \ref{thm2}. The following proposition is the main ingredient of the proof.
\begin{proposition}\label{prop}
Let $G$ be a $p$-group of nilpotency class 2 and $\Psi_2$ be the homomorphism given in the Proposition \ref{prop1}. Suppose $d\Big(\frac{G}{Z(G)}\Big) = \delta$ and $d(\gamma_2(G)) = \gamma$. Then, \[|\Ima(\Psi_2) | \geq p^{\sum\limits_{i= 2}^{min(\delta, \gamma+1)} (\delta-i)}.\]
\end{proposition}

$\allowbreak$
\begin{proof}
Choose $x_1, x_2, \ldots, x_{\delta} \in G$ such that 
\[\frac{G}{\Phi(G)Z(G)} = \gen{\overline{x_1}} \times \gen{\overline{x_2}} \times \cdots \times \gen{\overline{x_{\delta}}}.\] Let $U$ be the set $\{x_1, x_2, \ldots, x_{\delta}\}.$  We now choose a minimal generating set for $\frac{\gamma_2(G)}{\gamma_2(G)^p}$ in the following manner: Since the set $T = \{[x_i, x_j] \mid 1 \leq i < j \leq \delta\}$ generates $\gamma_2(G)$, we choose an element from this set, say $[x_{i_1^1}, x_{i_2^1}]$ such that $\overline{[x_{i_1^1}, x_{i_2^1}]} \neq 0 \in \frac{\gamma_2(G)}{\gamma_2(G)^p}$.
Define \[U_1 = \{x_{i_1^1}, x_{i_2^1}\},\] and 
\[V_1 = \Big\langle\overline{[x_{i_1^1}, x_{i_2^1}]}\Big\rangle \leq \frac{\gamma_2(G)}{\gamma_2(G)^p}.\] Suppose $U_j$ and $V_j \leq \frac{\gamma_2(G)}{\gamma_2(G)^p}$ have been defined. To define $U_{j+1}$ and $V_{j+1}$, we check if there exists any element $[y, z] \in T$, $y, z \in U$ such that $U_j \cap \{y, z\} \neq \phi$ and $V_j < \Big<V_j, \overline{[y,z]}\Big>$. If such an elements exists in $T$, say $[x_{i_1^{j+1}}, x_{i_2^{j+1}}]$, then we define 
\[U_{j+1} = U_j \cup \{x_{i_1^{j+1}}, x_{i_2^{j+1}}\},\] and 
\[V_{j+1} = \bigg<V_j, \overline{[x_{i_1^{j+1}}, x_{i_2^{j+1}}]}\bigg> \leq \frac{\gamma_2(G)}{\gamma_2(G)^p}.\] 
Otherwise we choose any element $[x_{i_1^{j+1}}, x_{i_2^{j+1}}] \in T$, $x_{i_1^{j+1}}, x_{i_2^{j+1}} \in U$ such that 
\[V_{j} < \Big<V_j, \overline{[x_{i_1^{j+1}}, x_{i_2^{j+1}}]}\Big>,\] and define 
\[U_{j+1} = \{x_{i_1^{j+1}}, x_{i_2^{j+1}}\},\] and 
\[V_{j+1} = \Big<V_j, \overline{[x_{i_1^{j+1}}, x_{i_2^{j+1}}]}\Big>.\] 
Clearly \[V_{\gamma} = \bigg<\overline{[x_{i_1^1}, x_{i_2^1}]}, \ldots, \overline{[x_{i_1^{\gamma}}, x_{i_2^{\gamma}}]}\bigg>  = \frac{\gamma_2(G)}{\gamma_2(G)^p}.\] 

Now suppose that for $j = k_1, k_2, \ldots, k_t$, $U_j \cap U_{j+1} = \phi$ 
and also that these are the only such numbers.
Denote the element \[\bar{x} \otimes \overline{[y,z]} + \bar{y} \otimes \overline{[z,x]} + \bar{z} \otimes \overline{[x,y]} \in \frac{G}{\Phi(G)Z(G)} \otimes \frac{\gamma_2(G)}{\gamma_2(G)^p}\] by $(x,y,z)$, and define 
\[W_{j} = \{(x, x_{i_1^{j}}, x_{i_2^{j}}) \mid x \in U \backslash U_{j} \}\] 
for $j = 1, \ldots, \gamma$.

We claim that $W_{1} \cup  W_{2} \cup \ldots \cup W_{\gamma}$ minimally generates $\Big<W_{1} \cup  W_{2} \cup \ldots \cup W_{\gamma}\Big>$. 
To see this, note that  $\frac{G}{\Phi(G)Z(G)}$ and $\frac{\gamma_2(G)}{\gamma_2(G)^p}$ are elementary abelian $p$-groups, and therefore can be considered as vector spaces over field $\mathbb{Z}/p\mathbb{Z}$ with bases $\{\overline{x_1}, \ldots, \overline{x_{\delta}}\}$ and $\{\overline{[x_{i_1^1}, x_{i_2^1}]}, \ldots, \overline{[x_{i_1^{\gamma}}, x_{i_2^{\gamma}}]}\}$ respectively.  It follows that the set $\{\overline{x_i} \allowbreak \otimes \overline{[x_{i_1^{j}}, x_{i_2^{j}}]} \mid 1 \leq i \leq \delta, 1 \leq j \leq \gamma\}$ form a basis for the tensor product $\frac{G}{\phi(G)Z(G)} \otimes \frac{\gamma_2(G)}{\gamma_2(G)^p}$. 
Now take an element $(x, x_{i_1^1}, x_{i_2^1}) \in W_1$. 
The presence of the term $\overline{x} \otimes \overline{[x_{i_1^1}, x_{i_2^1}]}$ in the expression $(x, x_{i_1^1}, x_{i_2^1})$ ensures that $(x, x_{i_1^1}, x_{i_2^1}) \notin \Big<W_1 \backslash (x, x_{i_1^1}, x_{i_2^1})\Big>$. This shows that $W_1$ minimally generates $\gen{W_1}$. If $k_1 > 1$, suppose for $j \leq k_1-1$, $W_1 \cup W_2 \ldots \cup W_j$ minimally generates $\Big<W_1 \cup W_2 \ldots \cup W_j\Big>$. Take an element $(x, x_{i_1^{j+1}}, x_{i_2^{j+1}}) \in W_{j+1}$. By definition $x \in U \backslash U_{j+1}$. Therefore if
\[(x, x_{i_1^{j+1}}, x_{i_2^{j+1}}) \in \bigg<W_1, W_2, \ldots, W_j, W_{j+1} \backslash (x, x_{i_1^{j+1}}, x_{i_2^{j+1}})\bigg>,\]
then 
\[\overline{x} \otimes \overline{[x_{i_1^{j+1}}, x_{i_2^{j+1}}]} = \overline{x} \otimes \overline{[x_{i_1^{1}}, x_{i_2^{1}}]^{\alpha_1} \cdots [x_{i_1^{j}}, x_{i_2^{j}}]^{\alpha_j}}.\] 
Since $\overline{x} \neq 0 \in \frac{G}{\Phi(G)Z(G)}$, we get that \[ \overline{[x_{i_1^{j+1}}, x_{i_2^{j+1}}]} = \overline{[x_{i_1^{1}}, x_{i_2^{1}}]^{\alpha_1} \cdots [x_{i_1^{j}}, x_{i_2^{j}}]^{\alpha_j}}.\] 
But this is not possible because $V_j < \Big<V_j, \overline{[x_{i_1^{j+1}}, x_{i_2^{j+1}}]}\Big>.$ This shows that $W_{1} \cup  W_{2} \cup \ldots \cup W_{k_1}$ minimally generates $\gen{W_{1} \cup  W_{2} \cup \ldots \cup W_{k_1}}$. 

Now if $\gamma > 1$, suppose that $k_{t+1} = \gamma$, and also that $W_{1} \cup  W_{2} \cup \ldots \cup W_{k_j}$ for $1 \leq j \leq t$ minimally generates $\Big<W_{1} \cup  W_{2} \cup \ldots \cup W_{k_j}\Big>$. Let $(x, x_{i_1^{k_j+1}}, x_{i_2^{k_j+1}}) \in W_{k_j+1}$. By the definition of  $W_{k_j+1}$, we have $x \in U \backslash U_{k_j+1} \ (i.e. \ U \backslash \{x_{i_1^{k_j+1}}, x_{i_2^{k_j+1}}\})$. First assume that $x \in U_{k_1} \cup \cdots \cup U_{k_j}$ and let $l_1, l_2, \ldots, l_r$ be all those numbers less than $k_j+1$ such that $ x = x_{i_{s_1}^{l_1}} = x_{i_{s_2}^{l_2}} = \cdots = x_{i_{s_r}^{l_r}}$ where $s_i's$ are either 1 or 2. Without loss of generality, we can assume that  $s_i = 1$ for $i = 1, 2, \ldots, r$. 
Now if
\[(x, x_{i_1^{k_j+1}}, x_{i_2^{k_j+1}}) \in \bigg<W_1, W_2, \ldots  W_{k_j+1}\backslash (x, x_{i_1^{k_j+1}}, x_{i_2^{k_j+1}})\bigg>,\]
 then 
 \begin{eqnarray*}
 \overline{x} \otimes \overline{[x_{i_1^{k_j+1}}, x_{i_2^{k_j+1}}]} & = &
  \overline{x} \otimes \overline{[x_{i_1^{1}}, x_{i_2^{1}}]^{\alpha_1} \cdots [x_{i_1^{k_j}}, x_{i_2^{k_j}}]^{\alpha_{k_j}}} \overline{[x_{i_2^{l_1}}, y_1]^{\beta_1}  
  \cdots [x_{i_2^{l_r}}, y_r]^{\beta_r}}
  \end{eqnarray*}
for some $y_1, y_2, \ldots, y_r \in G$ and for some $\alpha_i, \beta_k \in \mathbb{Z}$. Since $\overline{x} \neq 0 \in \frac{G}{\Phi(G)Z(G)}$ we have 
\begin{eqnarray*}
 \overline{[x_{i_1^{k_j+1}}, x_{i_2^{k_j+1}}]} & = & \overline{[x_{i_1^{1}}, x_{i_2^{1}}]^{\alpha_1} \cdots [x_{i_1^{k_j}}, x_{i_2^{k_j}}]^{\alpha_{k_j}} [x_{i_2^{l_1}}, y_1]^{\beta_1} [x_{i_2^{l_2}}, y_2]^{\beta_2} \cdots [x_{i_2^{l_r}}, y_r]^{\beta_r}}.
\end{eqnarray*}
But then, for some $y \in U$ and for some $i$ , 
\[V_{k_j} < \Big<V_{k_j}, \overline{[x_{i_2^{l_i}}, y]}\Big>.\] This contradicts the way we have chosen the basis for $\frac{\gamma_2(G)}{\gamma_2(G)^p}$. Therefore, we can now assume that that $x \in U \backslash (U_{k_1} \cup \cdots \cup U_{k_j} \cup U_{k_j+1})$. Next, if 
\[(x, x_{i_1^{k_j+1}}, x_{i_2^{k_j+1}}) \in \Big<W_1, W_2, \ldots  W_{k_j+1}\backslash (x, x_{i_1^{k_j+1}}, x_{i_2^{k_j+1}})\Big>,\] 
we get that 
\begin{eqnarray*}
\overline{x} \otimes \overline{[x_{i_1^{k_j+1}}, x_{i_2^{k_j+1}}]} & = & \overline{x} \otimes \overline{[x_{i_1^{1}}, x_{i_2^{1}}]^{\alpha_1} \cdots  [x_{i_1^{k_j}}, x_{i_2^{k_j}}]^{\alpha_{k_j}}}
\end{eqnarray*}
for some $\alpha_i \ (1 \leq i \leq k_j) \in \mathbb{Z}$, which is again not possible. This shows that $W_{1} \cup  W_{2} \cup \ldots \cup W_{k_j+1}$  minimally generates $\Big<W_{1} \cup  W_{2} \cup \ldots \cup W_{k_j+1}\Big>$. If $k_{j+1} > k_j+1$, suppose that $W_{1} \cup  W_{2} \cup \ldots \cup W_{k_j+i}$ for $i \leq k_{j+1}-k_j-1$ minimally generates $\Big<W_{1} \cup  W_{2} \cup \ldots \cup W_{k_j+i}\Big>$. With the same idea applied in the above argument, it can be shown that  $W_{1} \cup  W_{2} \cup \ldots \cup W_{k_j+i+1}$ minimally generates $\Big<W_{1} \cup  W_{2} \cup \ldots \cup W_{k_j+i+1}\Big>$. This shows that $W_{1} \cup  W_{2} \cup \ldots \cup W_{k_{j+1}}$ minimally generates $\Big<W_{1} \cup  W_{2} \cup \ldots \cup W_{k_{j+1}}\Big>$. Therefore, we can now conclude that $W_{1} \cup  W_{2} \cup \ldots \cup W_{\gamma}$ minimally generates $\Big<W_{1} \cup  W_{2} \cup \ldots \cup W_{\gamma}\Big>$ and hence hence $W_{1} \cup  W_{2} \cup \ldots \cup W_{\gamma}$ is a linearly independent set. 

Now, since $|U_1| = |U_{k_i+1}| = 2$ for $1 \leq i \leq t$, it follows, by the definition of $W_1$ and $W_{k_i+1}$, that  $|W_1|$ = $|W_{k_i+1}| = \delta-2$. For $1 \leq j \leq k_i -1$, note that, if $|W_j| \geq m$ then $|W_{j+1}| \geq m-1$, because $U_j \cap \{x_{i_1^{j+1}}, x_{i_2^{j+1}}\} \neq \phi$. Therefore it easily follows that 
\[|W_{1} \cup  W_{2} \cup \ldots \cup W_{\gamma}| \geq \sum\limits_{i= 2}^{min(\delta, \gamma+1)} (\delta-i).\] 

Now, by the universal property of tensor products there exists a homomorphism $\eta$, such that the following diagram commutes. 

\begin{center}
                        $ \begin{CD}
\frac{G}{\gamma_2(G)Z(G)} \times \gamma_2(G)         @>\phi >>  \frac{G}{\gamma_2(G)Z(G)} \otimes \gamma_2(G)    \\
@V \mathcal{P} VV             @V \eta VV\\
\frac{G}{\Phi(G)Z(G)}\times \frac{\gamma_2(G)}{\gamma_2(G)^p}          @> \theta >>     \frac{G}{\Phi(G)Z(G)} \otimes \frac{\gamma_2(G)}{\gamma_2(G)^p},\\  
\end{CD}$
\end{center}

\vspace{.2cm}
where \[\mathcal{P}(\overline{x}, y) = (\overline{x}, \overline{y}) \ \ \text{for} \ x \in G, y \in \gamma_2(G),\]

\[\phi(\overline{x}, y) = \overline{x} \otimes y,\ \ \  \text{and} \ \ \ 
\theta(\overline{x}, \overline{y}) = \overline{x} \otimes \overline{y}.\] 

 \vspace{.1cm} 
 Therefore, we have
 \[\eta(\overline{x} \otimes y) =  \overline{x} \otimes \overline{y}.\]
  
 Since the preimage of $W_{1} \cup  W_{2} \cup \ldots \cup W_{\gamma}$ under $\eta$ is isomorphic to a subgroup of $\Ima(\Psi_2)$ we get that \[|\Ima(\Psi_2)| \geq p^{\sum\limits_{i= 2}^{min(\delta, \gamma+1)} (\delta-i)}.\]
 This completes the proof.
\end{proof} 



We are now ready to prove The Theorem \ref{thm2}.\\

\noindent \textbf{\textit{Proof of Theorem \ref{thm2}}}: Let $\Psi_2$ be the homomorphism  given in the Proposition \ref{prop1} and $\overline{\Psi}_2$ be the similarly defined homomorphism associated with the group $G/\gamma_3(G).$ Also, let $d(G/Z(G)) = \delta$. Since $G/\gamma_3(G)$ is a group of nilpotency class 2, we get from the Proposition \ref{prop} that 
\[|\Ima(\overline{\Psi}_2) | \geq p^{\sum\limits_{i= 2}^{min(\delta, \gamma+1)} (\delta-i)}.\]
It is easy to see that $|\Ima(\Psi_2)| \geq |\Ima(\overline{\Psi}_2)|.$ The Theorem \ref{thm2} now follows from Equation \ref{eq1}.\\

We now proceed towards the proof of the Corollary \ref{cor2}. The proof makes use of the following two results.

\begin{proposition}\label{prop3}
Let $G$ be a p-group ($p$ odd)  of nilpotency class 2 with $G/\gamma_2(G)$ elementary abelian, $d(G) =d$ and  $|G^p| =p^t$. Let $V$ be the subgroup of $\gamma_2(G) \otimes G/\gamma_2(G)$  generated by all elements of the form $x^p \otimes x\gamma_2(G)$ for $x \in G$. Then $|V| =p^{\frac{1}{2}t(2d-t+1)}.$

\end{proposition}

The proof of the proposition follows exactly along the same lines as \cite[proposition 3.3]{Rai2}.

\begin{theorem}\cite[Theorem 2.5.6]{Karp}\label{thma}
Let $Z$ be a central subgroup of a finite group $G$. Then there exists the following exact sequence
\[Z \otimes \frac{G}{\gamma_2(G)} \mapsto M(G) \mapsto M(G/Z) \mapsto \gamma_2(G) \cap Z \mapsto 1,\]
where the map $\alpha: Z \otimes \frac{G}{\gamma_2(G)} \mapsto M(G)$ is defined as follows:
Let $G$ be given by $F/R$ for some free group $F$ and its normal subgroup $R$, and $Z$ be identified as $T/R$. After identifying $Z \otimes G$ and $M(G)$ as  $T/R \otimes F/\gamma_2(F)R$ and $\gamma_2(F) \cap R/[F,R]$ respectively, $\alpha$ is defined by 
\[\alpha(xR \otimes yR\gamma_2(F)) = [x,y][F,R].\]

\end{theorem}

We are now ready to prove Corollary \ref{cor2}.\\

\noindent \textbf{\textit{Proof of Corollary \ref{cor2}}}. For all finite $p$-groups, the lower bound is a well known fact [Corollary 3.2, \cite{Jones}]. The upper bound is a direct consequence of the Theorem \ref{thm2} and the fact that for the group $G$, $\gamma = k$ and $n=d+k$. 

Assume next that $G^p \cong \mathbb{Z}_p$. Since $\gamma_2(K) \leq K^2$ for any finite group $K$, for $p=2$ it follows that $G$ is either a  quaternion or a dihedral group of order 8. Therefore $M(G)$ is either trivial or of order 2. Hence we can assume now that $p$ is an odd prime. Consider the exact sequence from Theorem \ref{thma} for $Z = G^p$. We will show that the map  $\alpha:G^p \otimes \frac{G}{\gamma_2(G)} \mapsto M(G)$ is the trivial map in this case. To see this, let $G^p = \gen{g^p}$ for some $g \in G$. By the definition of $\alpha$ note that $x^p \otimes x\gamma_2(G) \in \kernal \ \alpha$ for all $x \in G$. Therefore  
$(gy)^p \otimes gy\gamma_2(G) \in \kernal \ \alpha$ for all $x, y \in G$. The bilinearity of the tensor product $\otimes$ implies that 
\[g^p \otimes y\gamma_2(G) + y^p \otimes g\gamma_2(G) \in \kernal \ \alpha.\]
But $y^p = (g^p)^m $ for some natural number $m$. Hence $g^p \otimes y\gamma_2(G) \in \kernal \ \alpha$ for all $y \in G$. As a result, $\alpha$ is the trivial map. It follows, from the exact sequence in \ref{thma}, that 
\[|M(G)| = \frac{|M(G/G^p)|}{|G^p|}.\]
Now apply the general bound obtained earlier in the corollary for the group $G/G^p$ to get the required bound in this case.

 Next suppose that $p$ is an odd prime and $G^p = \gamma_2(G)$. Applying the exact sequence in Theorem \ref{thma} again for $Z = G^p = \gamma_2(G)$, we get 
 \begin{equation}\label{eq2} |M(G)| \leq \frac{|M(G/\gamma_2(G)|}{|\gamma_2(G)|}\frac{|G^p \otimes G/\gamma_2(G)|}{|\kernal \ \alpha|}.\end{equation}
 Again by the definition of $\alpha$ it is evident that $x^p \otimes \overline{x} \in \kernal \ \alpha$ for all $x \in G.$ 
 Therefore from Proposition \ref{prop3} 
 \[|\kernal \ \alpha | \geq p^{\frac{1}{2}\gamma(2d-k+1)}.\] 
 Now putting  $|M(G/\gamma_2(G)| = p^{\frac{1}{2}d(d-1)}$, $|G^p \otimes G/\gamma_2(G)| = p^{dk}$ and the lower bound for $|\kernal \ \alpha|$ in the Inequality \ref{eq2} yield the desired result.\\

 We now prove Corollary \ref{cor3} which improves on Corollary \ref{cor2} in the case $G$ is a special $p$-group with $|G| \geq p^{13}$ and $|Z(G)| = p^3.$ The proof uses the well known connection between the capability of groups and the Schur multiplier.
 
 \vspace{.2cm}
 \textbf{\textit{Proof of Corollary \ref{cor3}}}: Let $Z^*(G)$ be the smallest central subgroup of $G$ such that $G/Z^*(G)$ is capable. Since $|G| \geq p^{13}$, by \cite[Theorem 5.7]{Mazur} $G$ is not capable. Therefore $Z^*(G)$ is non-trivial. Let $Z$ be a subgroup of $Z^*(G)$ of order $p$. Since $Z \leq Z(G) = \gamma_2(G)$, by \cite[Theorem 2.5.10]{Karp} we have
 \[|M(G)| = |M(G/Z)|/|Z|.\] Note that $G/Z$ is a nilpotent group of class 2 such that $G^p \leq \gamma_2(G)$ and $(G/Z)^p = \gamma_2(G/Z)$. This enables us to apply Corollary \ref{cor2} for the group $G/Z$ to obtain the desired result.

\vspace{.3cm}
We now prove Theorem \ref{thm3} which gives bounds on the Schur multiplier of $p$-groups of maximal class. These groups are an important class of finite $p$-groups and were first studied by Blackburn \cite{Blackburn0}. The following proof uses the well-known facts discovered by him.\\ 

\noindent \textbf{\textit{Proof of Theorem \ref{thm3}}} Let $P_1 = C_G(\gamma_2(G)/\gamma_4(G)).$ Choose arbitrary elements $s \in G \backslash P_1 \cup C_G(\gamma_{n-2}(G))$ and $s_1 \in P_1 \backslash \gamma_2(G)$. Then $s$ and $s_1$ generate $G$. If we define $s_i = [s_{i-1}, s]$ for $i \geq 2$, then $s_i \in \gamma_i(G) \backslash \gamma_{i+1}(G)$. Let $\Psi_i, i \geq 3$, be the map defined in the Proposition \ref{prop1}. Then
\begin{eqnarray*}
\Psi_i(\overline{s_1} \otimes \overline{s} \otimes \overline{s} \otimes \cdots \otimes \overline{s} \otimes \overline{s_1}) & = & \overline{[s_1, s, s, \cdots,s]_l}\overline{[s, s, \cdots, s, s_1]_r} \otimes \overline{s_1} + \overline{t} \otimes \overline{s}
\end{eqnarray*}
for some $t \in G$. \\
If $i$ is an odd integer, notice that 
\[\overline{[s_1, s, s, \cdots,s]_l} = \overline{[s, s, \cdots, s, s_1]_r}.\]
Since $p \neq 2$, it follows that $\Psi_i(\overline{s_1} \otimes \overline{s} \otimes \overline{s} \otimes \cdots \otimes \overline{s} \otimes \overline{s_1})$ is a non-identity element so that Im$\Psi_i$ is non-trivial. Using this fact, Equation \ref{eq1} yields the desired result.


\vspace{.3cm}



\begin{thebibliography}{999}



\bibitem{Berkovich}
Y. Berkovich, Z. Janko, Groups of Prime Power Order, Vol. 3 (de Gruyter, Berlin, 2011).

\bibitem{Blackburn0}
N. Blackburn, \emph{On a special class of p-groups}, Acta Math. 100 (1958), 49-92.

\bibitem{Blackburn}
N.  Blackburn, and  L.  Evens, \emph{Schur  multipliers  of p-groups},  J.  Reine  Angew.  Math. 309 (1979), 100-113.


\bibitem{Ellis1} 
G. Ellis and J. Wiegold, \emph{A bound on the Schur multiplier of a prime-power group},
Bull. Austral. Math. Soc. 60 (1999), 191-196.


\bibitem{Gaschutz}
W. Gaschutz, J. Neub\"{u}ser and Ti Yen, \emph{Uber den Multiplikator von p-Gruppen}, Math. Z. 100 (1967), 93-96.

\bibitem{Green}
J. A. Green, \emph{On the number of automorphisms of a finite p-group}, Proc. Roy. Soc. London 237 (1956), 574-581.
 
 \bibitem{Hatui}
 S. Hatui, \emph{Schur multipliers of special $p$-groups of rank 2}, J. Group Theory 23 (2020), 85–95.

\bibitem{Jones2}
M.R. Jones, \emph{Multiplicators of $p$-groups}, Math. Z. 127 (1972), 165-166.

\bibitem{Jones}
M.R. Jones, \emph{Some inequalities for the multiplicator of a finite group}, Proc. Amer. Math. Soc. 39 (1973), 450–456.

\bibitem{Jones1}
M.R. Jones \emph{Some inequalities for the multiplicator of a finite group II}, Proc. Amer. Math. Soc. 45, 167-172 (1974).

\bibitem{Karp}
G. Karpilovsky, \emph{The Schur multiplier}, London Math. Soc. Monogr, New Series no.2, (1987).

\bibitem{Khukhro}
E.I. Khukhro, \emph{$p$-Automorphisms of Finite $p$-Groups}, London Mathematical Society Lecture Note Series. 246, (1998).

\bibitem{Mazur}
M. Mazur, \emph{The epicenter of special $p$-groups of rank 2}, Unpublished work available at http://people.math.binghamton.edu/mazur/papers/pubp1.pdf

\bibitem{Moravec}
P. Moravec, \emph{On the Schur multipliers of finite p-groups of given coclass}, Israel J.Math. 185 (2011), 189-205.

\bibitem{Niroomand2}
P. Niroomand, \emph{On the order of Schur multiplier of non-abelian p-groups}, J. Algebra
322 (2009) 4479-4482.

\bibitem{Niroomand1}
P. Niroomand and F.G. Russo, \emph{An improvement of a bound of Green}, J. Algebra Appl., Vol. 11, No. 06, 1250116 (2012).


\bibitem{Rai1}
P.K. Rai, \emph{A note on the order of the Schur multiplier of $p$-groups}, International Journal of Algebra and Computation 27, (2017), 495-500.

\bibitem{Rai2}
P.K. Rai, \emph{ On the Schur multiplier of the Special $p$-groups}, Journal of Pure and Applied Algebra, 222, (2018), 316-322.

\bibitem{Vermani1}
L. R. Vermani, \emph{An exact sequence and a Theorem of Gaschutz, Neubuser and Yen on the
multiplicator},  J. London Math. Soc, 1 (1969), 95-100.

\bibitem{Vermani2}
L.R. Vermani, \emph{On the multiplicator of a finite group},  J. London Math. Soc, 1 (1974), 765-768.

\bibitem{Weigold1}
J. Wiegold, \emph{Multiplicators and groups with finite central factor-groups}, Math.
Z. 89 (1965), 345-347.

\bibitem{Weigold2}
J. Wiegold, \emph{Commutator subgroups of finite p-groups}, J. Austral. Math. Soc. 10, (1969), 480-484.


\end{thebibliography}
\end{document}